\font\eightrm=cmr8  
\font\eighttt=cmtt8
\magnification=\magstephalf

\parindent=0pt
\overfullrule=0in
 
\bf
\noindent
CHU'S 1303 IDENTITY IMPLIES BOMBIERI'S 1990 NORM-INEQUALITY
[Via An Identity of Beauzamy and D\'egot]
\medskip
\rm
(Appeared in the {\it Amer. Math. Monthly} {\bf 101}(1994), 894-896..)
\medskip
\it
\qquad\qquad
\qquad\qquad\qquad\qquad\qquad\qquad Doron ZEILBERGER\footnote{$^1$}
{\eightrm 
Department of Mathematics, Temple University,
Philadelphia, PA 19122, USA. 
{\eighttt zeilberg@math.temple.edu .} Supported in part by the NSF.
This note was written
while the author was on leave (Fall 1993) at the Institute for Advanced
Study, Princeton.  I would like to thank Don Knuth for a helpful suggestion.
}
\medskip
\qquad\qquad\qquad\qquad\qquad{\it Blessed are the meek: for they shall
inherit the earth (Matthew V.5)}
\medskip
\noindent
\medskip
\noindent
\rm
Inequalities are deep, while {\it equalities} are shallow. 
Nevertheless, it sometimes happens that a deep
inequality, {\bf A}, follows from a mere {\it equality} {\bf B}, which, in 
turn, follows from a more general, and {\it trivial}\footnote{$^2$}
{\eightrm Trivial to verify, not to conceive!} identity {\bf C}. 
\medskip
In this note we demonstrate this, following [3],
with {\bf A}$:=$ Bombieri's 
norm inequality[2]\footnote{$^3$}
{\eightrm It was needed by Beauzamy and Enflo in their research on
deep questions on Banach spaces. It also turned out to have far reaching
applications to computer algebra![1].
},
{\bf B}$:=$ an identity of Reznick[5], and {\bf C} $:=$ an identity of Beauzamy 
and D\'egot[3]. This exposition differs from the original only in the
punch line: I give a 1-line proof of {\bf C}, using Chu's identity.
 
\medskip
Let $P(x_1 , \dots , x_n)$ and $Q(x_1 , \dots , x_n)$ be two 
polynomials in $n$ variables:
 
$$
P= \sum_{ i_1 , \dots , i_n  \geq 0 } a_{ i_1 , \dots , i_n}
x_1^{i_1} \cdot \dots \cdot x_n^{i_n} \quad , \quad
Q= \sum_{i_1 , \dots , i_n  \geq 0 } b_{ i_1 , \dots , i_n}
x_1^{i_1} \cdot \dots \cdot x_n^{i_n}  \quad .
$$
 
The {\it Bombieri inner product}[2] is defined by
 
$$
[P,Q]:=
 \sum_{ i_1 , \dots , i_n  \geq 0 }( {i_1}! \dots {i_n}! ) \cdot
a_{ i_1 , \dots , i_n}
b_{ i_1 , \dots , i_n} \quad,
$$
 
and the {\it Bombieri norm}, by: $\Vert P \Vert := \sqrt{[P,P]} \quad .$
\medskip
{\bf Bombieri's Inequality A}: Let $P$ and $Q$ be any {\it homogeneous} 
polynomials in $(x_1 , \dots , x_n)$, then
 
$$
\Vert PQ \Vert \geq  \Vert P \Vert \Vert Q \Vert \quad .
$$
\medskip
In order to state {\bf B} and {\bf C}, we need to introduce the 
following notation.
$D_i := {{\partial} \over {\partial x_i}}$, ($i= 1, \dots ,n$),
$P^{( i_1 , \dots , i_n )} := D_1^{i_1} \dots D_n^{i_n} P$, and
for any polynomial $A(x_1, \dots , x_n)$, $A(D_1 , \dots , D_n)$
denotes the linear partial differential operator with constant
coefficients obtained by replacing $x_i$ by $D_i$.
\medskip
{\bf A} follows almost immediately from([5][3]):
\medskip
{\bf Reznick's Identity B:} For any polynomials $P$, $Q$ in $n$ variables:
 
$$
\Vert P Q \Vert ^2 =
 \sum_{i_1 , \dots , i_n \geq 0 } 
{{\Vert P^{(i_1 , \dots , i_n )} (D_1 , \dots , D_n) Q(x_1 , \dots , x_n ) 
\Vert^2} \over {i_1! \cdot \dots \cdot i_n!}} \quad.
$$
 
{\bf Beauzamy and D\'egot's Identity C:}
For any polynomials $P$,$Q$,$R$,$S$ in $n$ variables:
 
$$
[PQ,RS]=
\sum_{i_1 , \dots  , i_n  \geq 0} 
{{[R^{(i_1 , \dots , i_n)}(D_1 , \dots , D_n) Q(x_1 , \dots , x_n) ,
P^{(i_1 , \dots , i_n)}(D_1 , \dots , D_n) S(x_1 , \dots , x_n) ]}
\over
{( i_1! \dots i_n! )}} \quad.
$$
 
{\bf Proof of  B} $\Rightarrow$ {\bf A}:
Pick the terms for which  $ i_1 + \dots + i_n$ equals the (total)
degree of $P$, let's call it $p$, and note
that for those  $( i_1 , \dots , i_n )$,
$P^{( i_1 , \dots , i_n )} (x_1 , \dots , x_n ) =( i_1! \dots i_n!)
a_{i_1 , \dots , i_n}$, so
 
$$
\sum_{i_1 + \dots + i_n = p} 
{{\Vert P^{(i_1 , \dots , i_n )} (D_1 , \dots , D_n) Q(x_1 , \dots , x_n ) 
\Vert^2} \over {i_1! \cdot \dots \cdot i_n!}} 
=\sum_{i_1 + \dots + i_n = p} 
{{\Vert a_{i_1 , \dots , i_n} Q(x_1 , \dots , x_n ) 
\Vert^2} \cdot {( i_1! \cdot \dots \cdot i_n! ) }} 
$$
$$
=
\left ( \sum_{i_1 + \dots + i_n = p} 
{{ ( a_{i_1 , \dots , i_n})^2 } \cdot {( i_1! \cdot \dots \cdot i_n!) }} 
\right )
\Vert Q(x_1 , \dots , x_n ) 
\Vert^2 =  \Vert P \Vert^2 \Vert Q \Vert^2 \quad.
$$
\medskip
{\bf Proof of C} $\Rightarrow$ {\bf B:} Take $R=P$
and $S=Q$.
\medskip
{\bf Proof of C:} Both sides are
linear in $P$, in $Q$, in $R$, and in $S$, so it suffices to take
them all to be typical monomials,
($P=x_1^{p_1} \cdot \dots \cdot x_n^{p_n}$, and similarly for
$Q$,$R$, and $S$), for which
the assertion follows immediately  by applying Chu's[4]
identity\footnote{$^4$}
{\eightrm Rediscovered in the 18th century by Vandermonde. Proved by
counting, in two different ways, 
the number of ways of picking p lucky winners out of a set
of r boys and s girls.}
 
$$
\sum_{i \geq 0} { {r} \choose {i}} {{s} \choose {p-i}} = {{r+s} \choose {p}}
\quad ,
$$
 
to $r=r_t , s=s_t, p=p_t, (t= 1 \dots n)$, (using $i_t$ for $i$),
and taking their product. Q.E.D.
\medskip
\noindent
{\bf References}
\medskip
\noindent
1. B. Beauzamy, {\it Products of polynomials and a priori estimates
for coefficients in polynomial decompositions: A sharp result},
J. Symbolic Computation {\bf 13}(1992), 463-472.
 
2. B. Beauzamy, E. Bombieri, P. Enflo, and H.L. Montgomery,
{\it Products of polynomials in many variables}, J. Number Theory {\bf 36}
(1990), 219-245.
 
3. B. Beauzamy and  J. D\'egot, {\it Differential Identities}, 
I.C.M., Paris, preprint.
 
4. Chu Chi-kie, manuscript, 1303, China. 
(See J. Needham, Science and Civilization
in China, v. 3, Cambridge University Press, New York, 1959.)
 
5. B. Reznick, {\it An inequality for products of polynomials}, Proc.
Amer. Math. Soc. {\bf 117}(1993), 1063-1073.
\medskip
Nov. 3, 1993 ; Revised: April 5, 1994.
\bye